\documentclass[12pt]{article}
\usepackage{epsfig,amssymb,epic,eepic,color, graphics}
\usepackage[english]{babel}

\def\nd{\noindent}

\newenvironment{demo}{\nd {\bf Proof: }}{${}$\hfill $\diamond$ \medskip}

\newtheorem{theo}{Theorem}
\newtheorem{prop}{Proposition}
\newtheorem{defi}{Definition}
\newtheorem{lemm}{Lemma}
\newtheorem{coro}{Corollary}
\newtheorem{stat}{Statement}

\begin{document}
\sloppy
\date{}
\title{Dynamically ordered energy function for Morse-Smale diffeomorphisms on 3-manifolds}
\author{V.~Grines\thanks{N. Novgorod
State University, Gagarina 23, N. Novgorod,
603950 Russia, grines@vmk.unn.ru.}\and F.~
Laudenbach\thanks{Laboratoire de
math\'ematiques Jean Leray,  UMR 6629 du
CNRS, Facult\'e des Sciences et Techniques,
Universit\'e de Nantes, 2, rue de la
Houssini\`ere, F-44322 Nantes cedex 3,
France,
francois.laudenbach@univ-nantes.fr.}\and O.~
Pochinka\thanks{N. Novgorod State
University, Gagarina 23, N. Novgorod, 603950
Russia, olga-pochinka@yandex.ru.}}

\maketitle
\begin{abstract} This note deals with arbitrary  Morse-Smale diffeomorphisms in dimension 3 and  extends ideas from \cite{GrLaPo}, \cite{GrLaPo1}, where gradient-like case was considered.  We introduce a kind of Morse-Lyapunov function, called dynamically ordered, which fits well dynamics of diffeomorphism. The paper is devoted to finding
conditions to the existence of such 
an energy function, that is, a function whose  set of critical points coincides with the non-wandering set of the considered diffeomorphism. We show that the necessary and sufficient conditions to the existence of a dynamically ordered energy function reduces to the type of embedding of one-dimensional attractors and repellers, each of them is a union of zero- and one-dimensional unstable (stable) manifolds of periodic orbits of  a given Morse-Smale diffeomorphism on a closed 3-manifold. 
\end{abstract}

\section{Introduction and formulation of the results}

Let $M$ be a closed jrientable 3-manifold and $f:M\to M $ be a preserving orientation Morse-Smale diffeomorphism, that is: its nonwandering set $\Omega_f$ is finite, hence consists  of periodic points; $f$ is hyperbolic along $\Omega_f$ and the stable and unstable manifolds have transverse intersections. 

\begin{defi} A Morse function $\varphi:M\to \mathbb R$ is said to be a {\it Lyapunov  function} for $f$ if: 

1) $\varphi\bigl(f(x)\bigr)< \varphi(x)$ for every $x\not\in \Omega_f$; 

2) $\varphi\bigl(f(x)\bigr)=\varphi(x)$ for every $x\in  \Omega_f$.
\end{defi}

Sometimes we shall  speak of a Lyapunov function even when it is only defined on some domain
$N\subset M$, meaning that the above
conditions 1), 2) hold only for points $x\in
N$ such that $f(x)\in N$.

Let us recall that a $C^2$-smooth function
$\varphi:M\to\mathbb{R}$ is called a {\it
Morse function} if all its critical points
are non-degenerate. Using ideas from \cite{Sm74} it is possible to construct Lyapunov functions for $f$. For this aim, one considers the suspension of $f$, a 4-dimensional manifold which is fibered  over the circle and  is endowed with a Morse-Smale vector field $X$ transverse to the fibration. The method introduced by S. Smale in \cite{Sm74} for constructing Lyapunov function for Morse-Smale vector fields without closed orbits can be extended to suspension and allows one 
to construct a Lyapunov function $\Phi$ for $X$. The restriction of $\Phi$ to the base fibre, identified with $M$, is a Lyapunov function for diffeomorphism $f$. 

According to statement \ref{st0} below the periodic points of $f$ are critical points of its Lyapunov function $\varphi$  and the index of $\varphi$ at $p\in \Omega_f$ equals the dimension of $W^u_p$. At the same time any periodic point $p$ is a maximum of the restriction of $\varphi$ to the unstable manifold $W^u_p$ and a minimum of its restriction to the stable manifold $W^s_p$. If these extrema are non-degenerate then the invariant manifolds of $p$ are transversal to all regular level sets of $\varphi$ in some
neighborhood  of the point $p$. This local property is
useful for the construction of  a (global)
Lyapunov function. Next definition was introduced in \cite{GrLaPo1}.

\begin{defi} A Lyapunov function $\varphi:M\to\mathbb R$ for the Morse-Smale diffeomorphism $f:M\to M$ is called a {\it Morse-Lyapunov} function if every periodic point $p$ is a non-degenerate maximum (resp. minimum) of the restriction of $\varphi$  to the unstable (resp. stable) manifold $W^u_p$ (resp. $W^s_p$).
\end{defi}

Among the Lyapunov functions of $f$ those which are 
Morse-Lyapunov  form a generic set in the $C^\infty$-topology (see \cite{GrLaPo1}, theorem 1). 
In general, a Morse-Lyapunov function may have critical points which are not periodic points of $f$. 
 
\begin{defi} A Morse-Lyapunov function $\varphi$ is called an {\it energy function} for a Morse-Smale diffeomorphism $f$ if the set  of critical points of $\varphi$ coincides with $\Omega_f$. 
\end{defi} 
 
D. Pixton in  \cite{Pi1977} established an existence of energy function for any Morse-Smale diffeomorphisms given on closed smooth two-dimensional manifold and constructed  a gradient-like diffeomorphism on $\mathbb S^3$ which has no  energy function. According to S. Smale \cite{Sm74} any  Morse-Smale flow without closed trajectories (gradient-like flow)  given on closed smooth manifold of any dimension possesses by an  energy function. Thus there is an actual problem a finding of conditions to an existence of energy function for Morse-Smale diffeomorphisms. First step in this direction was made by the authors for gradient-like diffeomorphisms in the papers \cite{GrLaPo}, \cite{GrLaPo1}. 

\begin{figure}\begin{center}
\epsfig{file=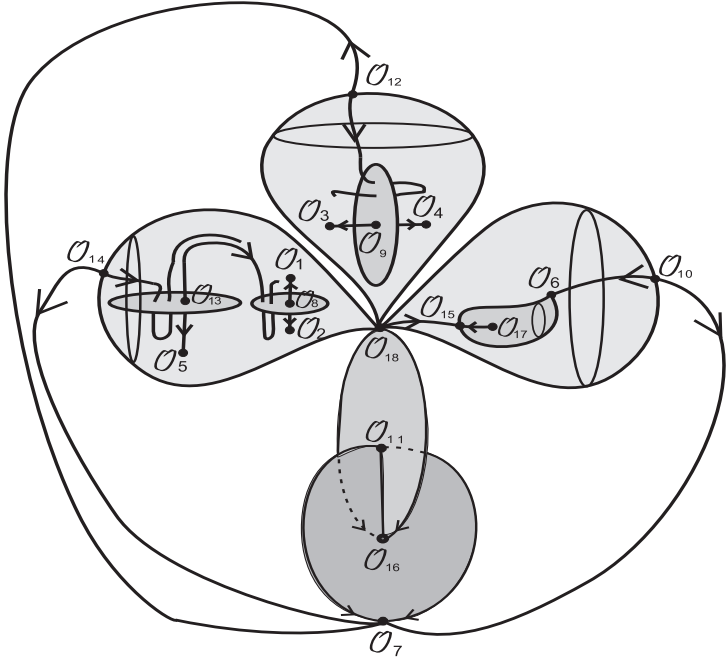, width=10 true cm, height=9. true cm}
\caption{Phase portrait of a Morse-Smale diffeomorphism $f:\mathbb S^3\to\mathbb S^3$ with dynamical numbering of the  periodic orbits}\label{beh}\end{center}
\end{figure}

Let us recall that a Morse-Smale diffeomorphism $f:M\to M$ is called {\it gradient-like} if for any pair of periodic points $x $, $y $ ($x\neq y $)
the condition $W^u_x\cap W^s_y\neq
\emptyset$ implies $\dim W^s_x<\dim W^s_y$. It follows from the definition  that a Morse-Smale diffeomorphism is gradient-like if and only if there are no {\it heteroclinic points} that  is,  intersection points of two-dimensional and one-dimensional invariant manifolds of different saddle  points.  Notice that two-dimensional invariant manifolds of different saddle points of a gradient-like diffeomorphism may have a non-empty intersection along the so-called {\it heteroclinic curves} (see figure \ref{beh}).

In \cite{GrLaPo}, \cite{GrLaPo1} (Theorem 4) we gave  necessary and sufficient conditions to the existence of a {\it self-indexing} energy function  for a Morse-Smale diffeomorphism $f:M\to M$ and showed that a non gradient-like diffeomorphisms do not possess a self-indexing energy function.  Here self-indexing means   $\varphi(p)=\dim W^u_p$ for every point $p\in \Omega_f$.

In the present paper we introduce the  notion of {\it dynamically ordered} Morse-Lyapunov function for  an arbitrary Morse-Smale diffeomorphism on 3-manifold. By using the above-mentioned arguments, such a function will  exist easily if it is not required to be an energy function.  We 
will show that the existence of such an energy function depends on how the one-dimensional attractors (and repellers)
embed into the ambient manifold.
More details are given below.

Let $f:M\to M$  be a Morse-Smale diffeomorphism. 
Following to S. Smale we introduce a partial order $\prec$ on the set 
of periodic orbits of $f$ in the following  way: 
$$\mathcal O_p\prec\mathcal O_r\iff W^s_{\mathcal O_p}\cap W^u_{\mathcal O_r}\neq\emptyset\,.$$
This definition means intuitively that all wandering points flow down along unstable manifolds to 
smaller elements. 
A sequence of different periodic orbits $\mathcal O_p=\mathcal O_{p_0},\mathcal O_{p_{1}},\dots,\mathcal O_{p_{k}}=\mathcal O_r$ ($k\geq 1$) such that $\mathcal O_{p_0}\prec \mathcal O_{p_{1}}\prec\dots\prec \mathcal O_{p_{k}}$ is called  a {\it chain of length $k$ connecting $\mathcal O_r$ to $\mathcal O_p$}. 
The maximum length of such chains is called, by J. Palis in \cite{Pa},  the {\it behaviour} of 
$\mathcal O_r$ relative to $\mathcal O_p$ and is denoted by $beh(\mathcal O_r\vert \mathcal O_p)$. For completeness it is assumed $beh(\mathcal O_r\vert \mathcal O_p)=0$ if $W^u_{\mathcal O_r}\cap W^s_{\mathcal O_p}=\emptyset$.

For each $q\in\{0,1,2,3\}$, denote $\Omega_q$ the subset of periodic points 
$r$ such that $\dim~W^u_{r}=q$ and denote $k_q$ the number of periodic orbits in the set $\Omega_q$. Set $k_f=k_0+k_1+k_2+k_3$ the number of all periodic orbits. For each periodic orbit $\mathcal O_r$ we set 
$q_{_{\mathcal O_r}}=\dim~W^u_{\mathcal O_r}$ and $b_{_{\mathcal O_r}}=\max\limits_{p\in\Omega_0}\{beh(\mathcal O_r\vert \mathcal O_p)\}$.

\begin{defi} \label{order} A numbering  of the periodic orbits: $\mathcal O_1,\dots,\mathcal O_{k_f}$ is called dynamical if it satisfies to following conditions: 

1) if $q_{\mathcal O_i}<q_{\mathcal O_j}$ then $i<j$;

2) if $q_{\mathcal O_i}=q_{\mathcal O_j}$ and $b_{\mathcal O_i}< b_{\mathcal O_j}$ then $i<j$.
\end{defi}

Notice that any dynamical numbering preserves the partial order $\prec$ (that is,$\mathcal O_i\prec \mathcal O_j$ implies $i\leq j$).  Indeed, as the  intersection  $W^s_{\mathcal O_i}\cap W^u_{\mathcal O_j}$ is transverse, the condition $\mathcal O_i\prec \mathcal O_j$ implies the inequality  $\dim~W^s_{\mathcal O_i}+\dim~ W^u_{\mathcal O_j}-3\geq 0$. Then $3-q_{\mathcal O_i}+q_{\mathcal O_j}-3\geq 0$ and, hence, $q_{\mathcal O_i}\leq q_{\mathcal O_j}$. If $q_{\mathcal O_i}<q_{\mathcal O_j}$ then $i<j$ due to 1). If $q_{\mathcal O_i}=q_{\mathcal O_j}$ then the condition $\mathcal O_i\prec \mathcal O_j$ implies or $\mathcal O_i=\mathcal O_j$ and, hence, $i=j$, either  $b_{\mathcal O_i}<b_{\mathcal O_{j}}$ and, hence, $i<j$ due to 2). 

On figure \ref{beh} it is represented a phase portrait of a Morse-Smale diffeomorphism $f:\mathbb S^3\to\mathbb S^3$ with $\Omega_f$ consisting of fixed points which are dynamically numerated.

\begin{defi} \label{dynod} Let $\mathcal O_1,\dots,\mathcal O_{k_f}$ be a dynamical
numbering of the  periodic orbits of $f$. A Morse-Lyapunov function $\varphi$ for $f$ is said to be dynamically ordered when $\varphi(\mathcal O_i)=i$ for $i\in\{1,\dots,k_f\}$.
\end{defi}

For each $i=1,\dots,k_1$, set   
$A_{i}=\bigcup\limits_{j=1}^iW^u_{\mathcal O_j}$.  It is known that  the set  $A_i$ is  an {\it attractor}, that is it has {\it a trapping neighborhood} $M_i$, which is a compact set such that $f(M_i)\subset int~M_i$ ($M_i$ is $f$-compressed) and $\bigcap\limits_{k\geq 0}f^k(M_i)=A_i$ (see, for example, \cite{Robinson-book99}). Denote by $r_i$ the number of saddles, by $s_i$ the number of sinks  and by $c_i$ the number of  connected components in  $A_i$. Set $g_i=c_i+r_i-s_{i}$.

Let us recall that a smooth compact orientable three-dimensional  manifold is called {\it  a handlebody of a genus $g\geq 0 $} if it is diffeomorphic to a manifold which is obtained from a closed 3-ball by an orientation reversing identification of $g$ pairs of pairwise disjoint closed 2-discs in  its  boundary. The boundary of such a handlebody is an orientable surface of genus $g$.

\begin{defi} \label{han} A trapping neighborhood $M_i$ of the attractor $A_i$ is called a handle if:
$M_i$ consists of  $c_i$ handlebodies. The sum  $g_{_{M_i}}$  of genera of all connected components of $M_i$ is called genus of the handle neighborhood.
\end{defi} 

Notice that for each $i=1,\dots,k_0$, the number $g_i$ equals  $0$, the attractor $A_i$ is zero-dimensional (as it consists of the sink orbits) and has a handle neighborhood $M_i$ of genus $g_i=0$ consisting of $c_i$ pairwise disjoint  3-balls (it follows, for example, from statement \ref{loc} below). For each $i=k_0+1,\dots,k_1$ the attractor $A_i$ contains an one-dimensional connected component, therefor we will say (taking liberty) that $A_i$ is one-dimensional attractor.

\begin{prop} \label{ttt} Each one-dimensional attractor $A_i$ of Morse-Smale diffeomorphism $f:M\to M$ has a handle trapping neighborhood $M_i$ with genus $g_{_{M_i}}\geq g_{i}$. 
\end{prop}

\begin{defi} \label{min} A handle neighborhood $M_i$ of one-dimensional attractor $A_i$ is said to be tight if:

1) $g_{_{M_i}}=g_{i}$;

2) $W^s_{\sigma}\cap M_i$ consists of exactly one two-dimensional closed disc for each saddle point $\sigma\in\mathcal O_i$.  

A one-dimensional attractor $A_i$ possessing tight trapping neighborhood $M_i$ is  said to be tightly embedded.
\end{defi}

By definition a {\it repeller}  for $f$ is an attractor for $f^{-1}$. Moreover, dynamical numbering of the orbits $\mathcal O_1,\dots,\mathcal O_{k_f}$ of a diffeomorphism  $f$ induces a dynamical numbering of the orbits $\tilde{\mathcal O}_1,\dots,\tilde{\mathcal O}_{k_f}$ of a diffeomorphism  $f^{-1}$ next way: $\tilde{\mathcal O}_i=\mathcal O_{k_f-i}$. Then  a one-dimensional repeller for $f$ is  said to be {\it tightly embedded} if it is such an attractor for $f^{-1}$ according to induced numbering. 

\begin{figure}\begin{center}
\epsfig{file=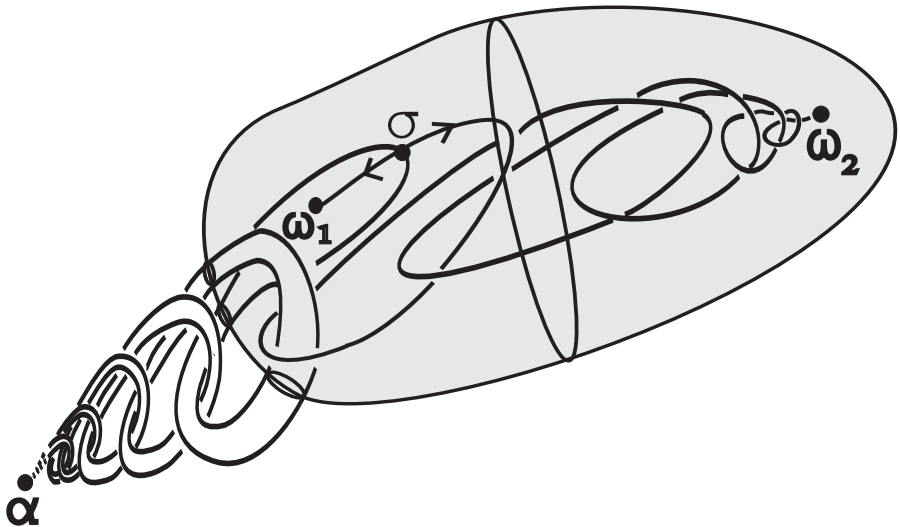, width=12 true cm, height=7 true cm}\caption{One-dimensional attractor of Pixton's  example is not tightly embedded} \label{ex}\end{center}
\end{figure}

Notice that the property for a one-dimensional attractor (repeller) to be tightly embedded gives a topological information about the  embedding of the unstable manifolds of its  saddle periodic points.  In the example which was constructed by D. Pixton in \cite{Pi1977}  the unique one-dimensional attractor $A_3=cl\,W^u_{\sigma}$ has the following property: $g_3=0$ but any 3-ball around $cl\,W^u_{\sigma}$ intersects $W^s_\sigma$ at more than one 2-disc (see figure \ref{ex}, where are drawn the phase portrait of Pixton's diffeomorpfism $f$ and a 3-ball). Hence,  this one-dimensional attractor is not tightly embedded.    

Our main results are the following theorems.

\begin{theo} \label{tt} If a  Morse-Smale diffeomorphism $f:M\to M$ possesses a dynamically ordered energy function, then all one-dimensional attractors and repellers of $f$ are tightly embedded. \end{theo}

\begin{defi} \label{strmin} A tight trapping neighborhood $M_i$ of  a one-dimensional attractor $A_i$ is called strongly tight if $M_i\setminus A_i$ is diffeomorphic to $\partial M_i\times(0,1]$.  A one-dimensional attractor $A_i$ possessing a strongly tight trapping neighborhood $M_i$ is 
said to be strongly tightly embedded. 
\end{defi}

\begin{theo} \label{iff} Let $f$ be a Morse-Smale diffeomorphism on a closed 3-manifold $M$. If all one-dimensional attractors and repellers of $f$ are strongly tightly embedded, then $f$ possesses a  dynamically ordered energy  function.
\end{theo}

Notice that the condition in  the last theorem is not necessary. For example in section 5 of paper \cite{GrLaPo1}  there was constructed a diffeomorphism on $\mathbb S^2\times\mathbb S^1$ possessing a  dynamically ordered energy  function, but whose  one-dimensional attractor and repeller  are not strongly tightly embedded.  

\begin{figure}\begin{center}
\epsfig{file=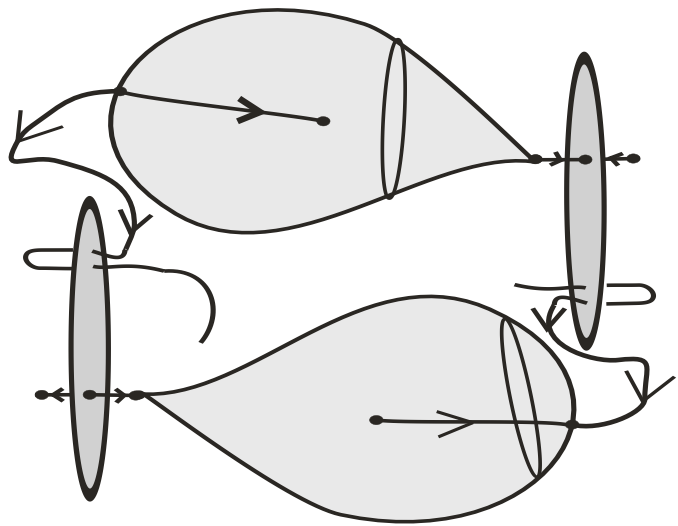, width=9 true cm, height=7 true cm}\caption{A Morse-Smale diffeomorphism without heteroclinic curves given on $\mathbb S^3$} \label{cha}\end{center}
\end{figure}

The next theorem states a criterion for the existence of some dynamically ordered energy  function for a  Morse-Smale diffeomorphism without heteroclinic curves given on $\mathbb S^3$. Methods 
from  \cite{BoGrPo2005} for realizing  Morse-Smale diffeomorphisms show that this class is not empty.
Moreover, it contains diffeomorphisms with chains of intersections of saddle invariant manifolds of arbitrary length (see figure \ref{cha}, where it is represented a phase portrait of a diffeomorphism from the class under consideration).  The criterion is based on paper \cite{BGMP}, where it is specified interrelation between topology of the 
ambient 3-manifold $M$ and structure of the non-wandering set of a Morse-Smale diffeomorphism without heteroclinic curves given on $M$. In particular, for any diffeomorphism of 
$\mathbb S^3$ without heteroclinic curves, the number $r$ of all saddles and the number $l$ of all sinks and sources satisfy the equality $r=l-2$ (see statement \ref{the} below). This equality implies that $g_i=0$ for any one-dimensional attractor $A_i$. Thus,  tightly embedded  attractor $A_i$ is strong tightly embedded. Applying results of theorems \ref{tt}, \ref{iff} we get next criterion. 

\begin{theo} \label{sph} A Morse-Smale diffeomorphism $f:\mathbb S^3\to\mathbb S^3$ without heteroclinic curves possesses a dynamically ordered energy function if and only if each one-dimensional attractor and repeller is tightly embedded. 
\end{theo}

\begin{center} ACKNOWLEDGMENTS
\end{center} 
V. Z. Grines and O. V. Pochinka acknowledge the
support of the grant of government of Russian Federation no. 11.G34.31.0039  for partial financial support. F. Laudenbach is supported by  the French program ANR ``Floer power".

\section{Auxiliary facts} 

In this section, we recall some statements that we need in the proof and   give references.

\begin{stat} ($\lambda$-lemma, \cite{Pa}). Let $p$ be a hyperbolic fixed point of a diffeomorphism  $f:M^n\to M^n$, $\dim W^u_{p}=\ell$, $0<\ell<n$, Let $B^u\subset W^u_{p}$ and  $B^s\subset W^s_{p}$ be small $\ell$-disc and  $(n-\ell)$-disc respectively
centered at  $p$. Let  $V:=B^u\times B^s$ be their product in a chart about $p$. Let $B$ be an $\ell$-disc transverse to $W^s_{p}$ at $x$. 
Then, for any $\varepsilon>0$, there exists  a positive integer $k_0$  such that the  connected component of $f^k(B)\cap V$ containing   $f^k(x)$ is 
$\varepsilon$-$C^1$-close to $B_u$ for each $k\geq k_0$. \label{lambda}
\end{stat}

\begin{defi} Let $F$ be a compact smooth surface properly embedded in a 3-manifold $W$ (that is,
 $\partial F\subset\partial W$). Then $F$ is called compressible in one from two following cases: 
 
1) there is a non contractible simple closed curve $c\subset int~F$ and smoothly embedded 2-disk  $D\subset int~W$ such that $D\cap F=\partial D=c$;

2) there is a 3-ball $B\subset int~W$ such that $F=\partial B$. 

The surface $F$ is said to be incompressible\footnote{It is well known to
topologists that a bicollared  surface, different from sphere, is incompressible if and only if  
the inclusion $S\hookrightarrow W$ induces an injection of  fundamental groups.} in $W$ if it is not compressible in $W$. 
\label{inncomp}
\end{defi}

\begin{stat} \label{wald32} {\rm (\cite{waldhausen}, corollary 3.2)} Let $S_g$ be an orientable surface of genus $g\geq 1$ and  let $F$ be an incompressible orientable surface properly embedded in $S_g\times[0,1]$ such that $\partial F\subset S_g\times\{1\}$. Then there is a 
surface $F_1\subset S_g\times\{1\}$ which is homeomorphic to $F$, such that 
$\partial F=\partial F_1$ and  $F\cup F_1$ bounds domain $\Delta$ in $S_g\times[0,1]$ such that  $cl\Delta$ is homeomorphic to $F\times[0,1]$, where $cl(\cdot)$ stands for the closure. 
\end{stat}

A particular case of statement \ref{wald32} is the following fact. 

\begin{coro} {\rm(\cite{GrMeZh}, theorem 3.3)} \label{the3}  Let $S_g$ be a closed orientable surface of genus $g\geq 1$ and let   surface $F\subset int(S_g\times[0,1])$ be a closed surface
which has genus $g$ and does not bound a domain in $S_g\times[0,1]$. 
Then $F$ is incompressible in $S_g\times[0,1]$ and the closure of each connected component of 
$S_g\times[0,1]\setminus F$ is homeomorphic to $S_g\times[0,1]$. 
\end{coro}

\begin{demo} According to the preceding statement, it is sufficient to check that $F$ is incompressible in $S_g\times[0,1]$. If $F$ is compressible, there exists some incompressible surface $F'$ whose genus $g'$ is less than $g$ and which still does not bound a domain in $S_g\times[0,1]$.  So $F'$ is not a sphere and $g'>0$. As $F'$ is incompressible, the preceding statement tells us that $F'$ is diffeomorphic to $S_g$. Contradiction.
\end{demo}

\begin{stat} {\rm (\cite{GrLaPo1}, lemma 3.3)} For any Morse-Smale diffeomorphism $f:M\to M$ 
we have that $1+|\Omega_1|-|\Omega_0|=1+|\Omega_2|-|\Omega_3|$,   where $|\cdot|$ stands for the cardinality. \label{g(f)}
\end{stat}

\begin{stat} \rm{(\cite{Mil1996}, theorem 5.2)} Let $M^n$ be a closed manifold, $\varphi:M^n\to\mathbb{R}$ be a Morse function, $C_q$ be  the number of all its critical points with index $q$, 
$\beta_q(M^n)$ be the $q$-th Betti number and $\chi(M^n)$ be  the Euler characteristic. Then 
$\beta_q(M^n)\leq C_q$ and $\chi(M^n)=\sum\limits_{q=0}^n(-1)^q C_q$.
\label{st10}
\end{stat}

\begin{stat} Let $\varphi:M^n\to\mathbb{R}$
be a Lyapunov function for a
Morse-Smale diffeomorphism $f:M^n\to
M^n$. Then

1) $-\varphi$ is Lyapunov function for
$f^{-1}$;

2) if $p$ is a periodic point of $f$
then   $\varphi(x)<\varphi(p)$ for every
$x\in W^u_p\setminus p$ and
$\varphi(x)>\varphi(p)$ for every
$x\in W^s_p\setminus p$;

3) if $p$ is a periodic point of $f$
then $p$ is a critical point of
$\varphi$ whose index is $\dim W^u_p$. \label{st0}
\end{stat}

\begin{stat} {\rm (\cite{GrLaPo1}, lemma 2.2)}
Let $f:M^n\to M^n$ be a Morse-Smale diffeomorphism 
on an $n$-dimensional manifold and let 
$\mathcal O$ be a periodic orbit.
For $p\in \mathcal O$, set $q=\dim W^u_p$. Then, there is some 
neighborhood $U$ and an energy function $\varphi:U\to\mathbb R$ for $f$ 
such that
 $(W^u_p\cap U)\subset Ox_1\dots x_q,~ (W^s_p\cap U)\subset Ox_{q+1}\dots x_n$
 for Morse coordinates 
$x_1,\dots,x_n$ of $\varphi$ near $p$.
\label{loc}
\end{stat}

\begin{stat} \label{the} {\rm (\cite{BGMP}, theorem)} Let $M$ be a three-dimensional closed, connected, orientable manifold. Let  $f:M\to M$ be any
Morse-Smale diffeomorphism without heteroclinic curves whose non-wandering set consists of $r$ saddles
and $l$ nodes (sinks and sources). Then $m=\frac{r-l+2}{2}$ is non negative  integer and  following facts hold: 

1) if  $m=0$, then   $M$ is the 3-sphere; 

2) if $m>0$, then   $M$ is the  connected sum
 of $m$ copies $\mathbb S^2\times\mathbb S^1$. 

\nd Conversely,   for any non negative  integers  $r, l, m$  such that $m=\frac{r-l+2}{2}$ is non negative  integer, there exists 3-manifold $M$ and some Morse-Smale diffeomorphism $f:M\to M$ with following properties: 

a) $M$ is 3-sphere if $m=0$ and $M$ is the  connected sum of $m$ copies $\mathbb S^2\times\mathbb S^1$ if $m>0$;

b) the  non-wandering set  of $f$ consists of $r$ saddles  and $l$ sinks and sources, the  wandering set of $f$ has no heteroclinic curves. 
\end{stat} 

\section{On one-dimensional attractors}

{\bf Proof of proposition \ref{ttt}} 

For each $i=k_0+1,\dots,k_1$ let us prove the  existence of handle neighborhood $M_i$ for one-dimensional attractor $A_i$ with $g_{_{M_i}}\geq g_i$.

\begin{demo} According to statement \ref{loc}, there is a neighborhood $U_{A_{k_0}}\subset  W^s_{A_{k_0}}$ of zero-dimensional attractor $A_{k_0}$ and an energy function  $\varphi_{A_{k_0}}:U_{A_{k_0}}\to\mathbb R$ for $f$ such that $\varphi_{A_{k_0}}(A_{k_0})=0$ and for small $\varepsilon>0$ each connected component of set $M_{k_0}=\varphi^{-1}_{A_{k_0}}((-\infty,\varepsilon])$ reads  $\{(x_1,x_2,x_3)\in U_{A_{k_0}} 
~:~ x^2_1+x^2_2+x^2_3\leq\varepsilon\}$ 
in local coordinates $x_1,x_2,x_3$. Then $M_{k_0}$ is trapping neighborhood of zero-dimensional attractor $A_{k_0}$, which is a union of $c_{k_0}$ pairwise disjoint 3-balls. By induction on $i=k_0+1,\dots,k_1$ we construct a handle trapping neighborhood $M_{i}$ for $A_{i}$.

Let $i=k_0+1$. Set $S_{k_0}=\partial M_{k_0}$. Without loss of generality we can suppose that $S_{k_0}$ intersects $W^u_{\mathcal O_{k_0+1}}$ transversely; let $n_{k_0}$  be the number of intersection points. Set $V_{k_0}=W^s_{\Omega_f\cap A_{k_0}}\setminus  A_{k_0}.$ Then the quotient $\hat V_{k_0}=V_{k_0}/f$ is made of the cobordism $M_{k_0}\setminus int\,f(M_{k_0})$
by gluing its boundaries by $f$. Hence, $\hat V_{k_0}$ is smooth orientable 3-manifold without boundary and natural projection $p_{_{k_0}}:V_{k_0}\to\hat V_{k_0}$ is cover. Then 
$p_{_{k_0}}(W^u_{\mathcal O_{k_0+1}})$ is a pair of knots which  intersects $p_{_{k_0}}(S_{k_0})$ transversely at $n_{k_0}$ points. Thus, there is a tubular neighborhood 
$\hat T^u_{k_0+1}\subset\hat V_{k_0}$ of  $p_{_{k_0}}(W^u_{\mathcal O_{k_0+1}})$ such that  $\hat T^u_{k_0+1}\cap p_{_{k_0}}(S_{k_0})$ consists of $n_{k_0}$ 2-discs.  

\begin{figure}\begin{center}
\epsfig{file=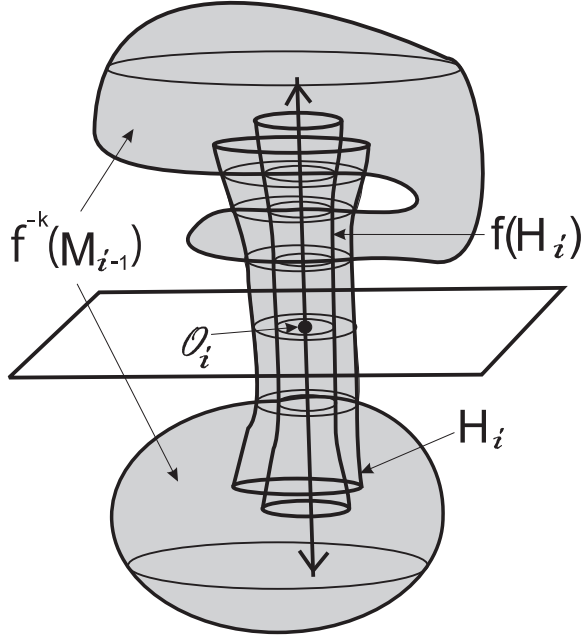, width=8 true cm, height=8.5 true cm}
\caption{Construction of a handle neighborhood for  one-dimensional attractor}\label{suflu2}\end{center}
\end{figure}

Set $T^u_{k_0+1}=p^{-1}_{k_0}(\hat T^u_{k_0+1})$. According to the $\lambda$-lemma (see statement \ref{lambda}), $T^u_{k_0+1}\cup W^s_{\mathcal O_{k_0+1}}$ is a neighborhood of $\mathcal O_{k_0+1}$. According to statement \ref{loc}, there are some neighborhood $U_{\mathcal O_{k_0+1}}\subset (T^u_{k_0+1}\cup W^s_{\mathcal O_{k_0+1}})$ of $ \mathcal O_{k_0+1}$ 
and an energy function $\varphi_{\mathcal  O_{k_0+1}}:U_{\mathcal O_{k_0+1}}\to\mathbb R$ for $f$ with $\varphi_{\mathcal O_{k_0+1}}(\mathcal O_{k_0+1})=0$. When $\varepsilon>0$  is small enough, each connected component of $H_{k_0+1}=\varphi^{-1}_{ \mathcal O_{k_0+1}}((-\infty,\varepsilon])$ reads $\{(x_1,x_2,x_3)\in U_{\mathcal O_{k_0+1}} 
~:~ -x^2_1+x^2_2+x^2_3\leq\varepsilon\}$ 
in local coordinates $x_1,x_2,x_3$. According to the $\lambda$-lemma, when $k\in\mathbb N$ is  large enough, $f^{-k}(S_{k_0})$  intersects both $H_{k_0+1}$ and $f(H_{k_0+1})$, its intersection with these domains consists of $n_{k_0}$ 2-discs and
$f(H_{k_0+1})\setminus int\,f^{-k}(M_{k_0})
\subset int\,H_{i}$ (see figure \ref{suflu2}). 
Thus, $M_{k_0+1}=f^{-k}(M_{k_0})\cup H_{k_0+1}$  
is a union of handlebodies, as it is obtained from the union of 3-balls $f^{-k}(M_{k_0})$ by gluing one-handles $H_{k_0+1}\setminus int\,f^{-k}(M_{k_0})$. Let us show that  $f(M_{k_0+1})\subset int\,M_{k_0+1}$. 

Indeed, it is true for a point $x\in f^{-k}(M_{k_0})$ as $f^{-k}(M_{k_0})$ is $f$-compressed and it is true for a point $x\in (H_{k_0+1}\setminus f^{-k}(M_{k_0}))$ as $f(H_{k_0+1})\setminus int\,f^{-k}(M_{k_0})\subset int\,H_{k_0+1}$.

Let us prove the equality $\bigcap\limits_{k\geq 0}f^k(M_{k_0+1})=A_{k_0+1}$. As $A_{k_0+1}\subset M_{k_0+1}$ and  $f^k(A_{k_0+1})=A_{k_0+1}$ for $k\in\mathbb Z$ then $A_{k_0+1}\subset\bigcap\limits_{k\geq 0} f^k(M_{k_0+1})$. Let us set $A'_{k_0+1}=\bigcap\limits_{k\geq 0} f^k(M_{k_0+1})$  and show that $A'_{k_0+1}=A_{k_0+1}$. Assume contrary: there is a point $x\in(A'_{k_0+1}\setminus A_{k_0+1})$. Due to theorem 2.3 in \cite{S3} there is a point  $p\in(\Omega_f\setminus A_{k_0+1})$ such that $x\in W^u_p$. As the set $A'_{k_0+1}$ is closed and invariant then $cl~(\mathcal O_x)\subset A'_{k_0+1}$ and, hence, $p\in A'_{k_0+1}$. This is a contradiction with the fact $ A'_{k_0+1}\subset W^s_{A_{k_0+1}\cap\Omega_f}$. 

Recall that we denote by $r_i$ the number of saddles, by $s_i$ the number of sinks, by $c_i$ the number of connected components of the attractor $A_i$ and set $g_i=c_i+r_i-s_{i}$. By the construction $M_{k_0+1}$ consists of $c_{k_0+1}$ 3-balls with 1-handles\footnote{Recall that a 3-dimensional 1-handle is the product of an interval with a 2-disc.}. Denote by $g_{_{M_{k_0+1}}}$ the sum of genera of connected components of $M_{k_0+1}$. Let us show that $g_{_{M_{k_0+1}}}\geq g_{k_0+1}$.  

In denotation above, the number of points in the orbit $\mathcal O_{{k_0+1}}$ equals $r_{{k_0+1}}- r_{k_0}$.  As $A_{{k_0+1}}= A_{k_0}\cup W^u_{\mathcal O_{{k_0+1}}}$ and $cl~W^u_{\mathcal O_{{k_0+1}}}\setminus W^u_{\mathcal O_{{k_0+1}}}\subset A_{k_0}$ then $c_{{k_0+1}}\leq c_{k_0}$. Denote by $l_{k_0+1}$ the number of connected components of the set $H_{{k_0+1}}\setminus int\,f^{-k}(M_{k_0})$. By the construction each of them is 1-handle and removing of $(l_{{k_0+1}}-(c_{k_0}-c_{{k_0+1}}))$ 1-handles from $M_{k_0+1}$ gives the set with the same $c_{k_0+1}$ connected components. Then the sum of genera $g_{_{M_{k_0}}}$ of $M_{k_0}$ can be calculate by formula $g_{_{M_{k_0}}}=g_{_{M_{k_0+1}}}-(l_{{k_0+1}}-(c_{k_0}-c_{{k_0+1}}))$.  
By the construction $l_{k_0+1}\geq(r_{{k_0+1}}- r_{k_0})$, hence $g_{_{M_{k_0}}}\leq g_{_{M_{k_0+1}}}-(r_{{k_0+1}}-r_{k_0}-(c_{k_0}-c_{{k_0+1}}))$ and $g_{_{M_{k_0+1}}}\geq g_{_{M_{k_0}}}+r_{{k_0+1}}-r_{k_0}-c_{k_0}+ c_{{k_0+1}}$. As $g_{_{M_{k_0}}}= g_{k_0}$ then $g_{_{M_{k_0+1}}}\geq c_{k_0}+r_{k_0}-s_{k_0}+r_{{k_0+1}}-r_{k_0}-c_{k_0}+ c_{{k_0+1}}=c_{{k_0+1}}+ r_{{k_0+1}}-s_{k_0}$. As $s_{k_0}=s_{k_0+1}$ then $g_{_{M_{k_0+1}}}\geq g_{k_0+1}$.

A smoothing of the set $M_{k_0+1}$ is  the required handle trapping neighborhood. 

Assuming that handle neighborhood for attractor $A_{i-1}$ already constructed, repeating construction above (changing $k_0$ by $i-1$), we construct $f$-compressed set $M_{i}=f^{-k}(M_{i-1})\cup H_{i}$, 
being a union of handle neighborhood $f^{-k}(M_{i-1})$ with 1-handles $H_{i}\setminus int\,f^{-k}(M_{i-1})$. It is similar proved that $M_{i}$ is required handle neighborhood.
\end{demo}

\vfill\eject

\begin{prop} \label{coon} The one-dimensional attractor  $A_{k_1}$ is connected.
\end{prop}
\begin{demo} Firstly, let us prove that any trapping neighborhood $M_{k_1}$ of $A_{k_1}$ is connected. Let us assume the contrary: $M_{k_1}$ is a union of pairwise disjoint closed sets $B_1$ and $B_2$.  As $M_{k_1}$ is $f$-compressed then without loss of generality we can suppose that $f(B_i)\subset int\,B_i,~i=1,2$. By  construction, $U_1=\bigcup\limits_{k>0}f^{-k}(int\,B_1)$, $U_2=\bigcup\limits_{k>0}f^{-k}(int\,B_2)$ are pairwise disjoint open sets
and $U_1\cup U_2=W^s_{\Omega_0\cup\Omega_1}$. On the other hand
$W^s_{\Omega_0\cup\Omega_1}=M\setminus W^s_{\Omega_2\cup\Omega_3}$ and, hence, $W^s_{\Omega_0\cup\Omega_1}$ is connected as  $\dim\,M=3$ and  $\dim\,W^s_{\Omega_2\cup\Omega_3}\leq 1$. This is a contradiction. 

Thus $A_{i}$ is connected as intersection of 
nested connected compact sets 
$M_{i}\supset f(M_{i})\supset\dots\supset f^k(M_{i})\supset\dots$. 
\end{demo}

\section{Necessary condition for existence of dynamically ordered energy function} \label{nas}

{\bf Proof of theorem \ref{tt}} 

Let us prove that if a Morse-Smale  diffeomorphism $f:M\to M$  has  a dynamically ordered energy function, then its one-dimensional 
attractors and repellers  are tightly embedded.

\begin{demo} Notice  that  $f$ and $f^{-1}$ possess  dynamically ordered energy functions simultaneously. Indeed, if $\varphi:M\to\mathbb R$ is such function for $f$ then $-\varphi:M\to\mathbb R$ is an energy function for $f^{-1}$ (see statement \ref{st0}) and $\tilde\varphi=k_f+1-\varphi:M\to\mathbb R$ is dynamically ordered energy function for $f^{-1}$. Therefore, it is enough to prove the fact for attractors.  

Let $\varphi:M\to\mathbb R$ be a dynamically ordered energy  function for $f:M\to M$, $i=k_0+1,\dots,k_1$ and $M_i=\varphi^{-1}([1,i+\varepsilon_{i}]),\varepsilon_{i}>0$. It follows from properties of dynamically ordered energy  function and statement  \ref{st0} that any orbit $\mathcal O_j$ with number $j\leq i$ belongs to $M_i$. Due to statement  \ref{st0}, $W^u_{\mathcal O_j}\subset M_i$. Thus  
$A_i\subset M_i$. It follows from definition of Lyapunov function that $f(M_i)\subset int~M_i$.  
Similar to proposition \ref{ttt} it is proved equality $\bigcap\limits_{k\geq 0}f^k(M_{i})=A_{i}$. Thus, $M_i$ is trapping neighborhood of attractor  $A_i$. Then $M_i$ has the same number of connected components as  $A_i$. Let us prove that there is $\varepsilon_{i}>0$ such that  $M_{i}$ is a  tight  neighborhood of $A_i$. 

As $\varphi$ is a Morse-Lyapunov function then there is $\varepsilon_{i}>0$ such that $W^s_{\sigma}\cap M_{i}$ consists of exactly one closed 2-disc for each saddle point $\sigma\in\mathcal O_{i}$. It follows from properties of dynamically ordered energy  function and statement  \ref{st0} that  $\varphi|_{M_i}$ has exactly $r_i+s_i$ critical points, among of them $s_i$ points have index $0$ and $r_i$ points  have index $1$. According to Morse theory,  $M_i$ is a union of of $s_{i}$ 3-balls with gluing of $r_i$ 1-handles and hence is a union of $c_i$ handlebodies. Denote by $g_{_{M_{i}}}$ the sum of genus of handlebodies from $M_{i}$.  According to statement \ref{st10}, $\chi(M_{i})=s_{i}-r_{i}$. It follows from Morse theory that $M_{i}$ has the homotopy type of a cellular complex consisting of $s_{i}$ zero-dimensional and $r_i$ one-dimensional cells, then $-g_{_{M_{i}}}+c_{i}=m_{i}-r_{i}$ or $g_{_{M_{i}}}=g_{i}$. 
\end{demo}

\section{Construction of a dynamically ordered energy function for $f$}

Now $f$ is a Morse-Smale diffeomorphism on a closed 3-manifold $M$  and  its one-dimensional attractors and repellers are strongly tightly embedded. Construction of a dynamically ordered energy function for $f$ is based on technical lemmas of next section. 

Recall that, by assumption of theorem \ref{iff}, each one-dimensional attractor $A_i,~i=k_0+1,\dots,k_1$ is strongly tightly embedded and, hence, has a handle neighborhood $M_i$ of genus $g_i$ such that  $M_i\setminus A_i$ is homeomorphic to $S_i\times(0,1]$, where $S_i=\partial M_i$, and for each point $\sigma \in \mathcal O_i$ the intersection $W^s_{\sigma}\cap M_i$ consists of exactly one 2-disk. Set $D_i=M_i\cap W^s_{\mathcal O_i}$. 
According to statement \ref{loc}, for each zero-dimensional attractor $A_i,~i=1,\dots,k_0$, there is a handle neighborhood of genus $g_i=0$ which is a union of $c_i$ 3-balls, we will denote it by $M_i$ and set $S_i=\partial M_i$. 

For $i=1,\dots,k_1$ set $K_i=M_i\setminus int~f(M_i)$, $N_i=W^s_{A_i\cap\Omega_f}$ and $V_i=N_i\setminus A_i$. According to ring hypothesis and corollary \ref{the3}, $K_i$ is diffeomorphic to $S_i\times[0,1]$. As $V_i=\bigcup\limits_{n\in\mathbb Z}f^n(K_i)$ then $V_i$ is diffeomorphic to $S_i\times\mathbb R$.

\subsection{Extension of Lyapunov functions}

\begin{defi} Let $D$ be a  subset of $M$ which is diffeomorphic to product $S\times[0,1]$ for some (possibly non connected) surface  $S$. Then $D$  is said to be an $(f,S)$-compressed product  when there is a  diffeomorphism $g: D\to  S\times [0,1]$ such that $g^{-1}(S\times \{t\})$ bounds an 
$f$-compressed domain in $M$ for any $t\in  [0,1]$.    \end{defi}

\begin{prop} \label{f,S} Let $D$ be an $(f,S)$-compressed product. Then for any values $d_0<d_1$ there is an energy function $\varphi_{_D}:D\to\mathbb R$ for $f|_D$ such that  $\varphi_{_D}(g^{-1}(S\times \{0\}))=d_0$ and $\varphi_{_D}(g^{-1}(S\times \{1\}))=d_1$. 
\end{prop}
\begin{demo} The desired function  $\varphi_{_D}:D\to\mathbb R$ is defined by formula  $\varphi_{_D}(x)=d_0+t(d_1-d_0)$ for $x\in g^{-1}(S\times \{t\}),~t\in[0,1]$. 
\end{demo}

\begin{lemm} Let $i\in\{1,\dots,k_1\}$ and $P_i$, $Q_i$ be handle neighborhoods of genus $g_i$ of the attractor $A_i$. If there is a dynamically ordered energy function $\varphi_{_{Q_i}}:Q_i\to\mathbb R$ for $f$ with $S_{Q_i}=\partial Q_i$ as a level set then there is a dynamically ordered energy function $\varphi_{_{P_i}}:P_i\to\mathbb R$ for $f$ with $S_{P_i}=\partial P_i$ as a level set.
\label{ww}
\end{lemm}
\begin{demo} We follow to scheme of the proof of lemma 4.2 from \cite{GrLaPo}. Give some remarks. 

Without loss of generality we assume that $Q_i\subset int~P_i$ (in the opposite case, instead pair $(Q_i,\varphi_{_{Q_i}})$ we can use pair $(f^n(Q_i),\varphi_{_{f^n(Q_i)}})$, where $f^n(Q_i)\subset int~P_i$ and $\varphi_{_{f^n(Q_i)}}=\varphi_{_{Q_i}}f^{-n}$). 
As $V_i$ is diffeomorphic to $S_{i}\times\mathbb R$ then, according to ring hypothesis and corollary \ref{the3}, $G_i={P}_i\setminus int\,{Q}_i$ is a product. As handle neighborhoods $f^n(Q_i)$ and  $f^n(P_i)$ contain the attractor $A_i$ for each $n\in\mathbb Z$ then the surfaces $f^n(S_{Q_i})$ and  $f^{n}(S_{P_i})$ do not bound domains in  $V_i$ and, hence, are incompressible due to corollary  \ref{the3}. Now let us construct the function $\varphi_{_{P_i}}$, for this aim we consider  two cases: 1) $S_{P_i}\cap(\bigcup\limits_{n>0}f^{-n}
({S_{Q_i}}))=\emptyset$ and 2) $S_{P_i}\cap(\bigcup\limits_{n>0}f^{-n}
({S_{Q_i}}))\neq\emptyset$.

In case 1), let $m$ be the first positive integer such that  $f^m({P_i})\subset int\,{Q}_i$. If $m=1$, then $G_i$ is $(f,S_{i})$-compressed product and proposition \ref{f,S} yields the required function as extension of the function $\varphi_{_{Q_i}}$ to $G_i$.

If $m>1$, the surfaces $f({S}_{P_i}), f^{2}({S}_{P_i}), \dots,f^{m-1}({S}_{P_i})$  
are mutually ``parallel'', that is: two by two they bound a product cobordism (according to ring hypothesis and corollary \ref{the3}). Therefore they subdivide $G_i$ in $(f,S_{i})$-compressed products and, hence, proposition \ref{f,S} yields the required function as extension of the function $\varphi_{_{Q_i}}$ to $f^{m-1}(P_i)\setminus int~Q_i,f^{m-2}(P_i)\setminus int~f^{m-1}(P_i),\dots,P_i\setminus int~f(P_i)$ in series.  

In case 2),  without loss of generality we may assume that $S_{P_i}$ is transverse to $\bigcup\limits_{n>0}f^{-n}({S}_{Q_i})$,
which implies that there is a finite family $\mathcal C$ of  intersection curves. We are going to describe a process of decreasing the number of intersection curves by an isotopy of ${Q}_i$ among handle neighborhoods of genus $g_i$ possessing a dynamically ordered energy function for $f$ which is constant on the boundary of the neighborhood.

Firstly we consider all intersection curves from $\mathcal C$ which are homotopic to zero in $S_{P_i}$. Let $c$ be an  innermost such curve. Then there is a disc $\delta\subset S_{P_i}$ which is bounded by $c$ and such that $int\,\delta$ contains no curves from the family $\mathcal C$. As $c\subset f^{-n}(S_{Q_i})$ for some $n$ and $f^{-n}(S_{Q_i})$ is incompressible in $V_{i}$, then $c$ bounds a disc $d\subset f^{-n}(S_{Q_i})$. Then  2-sphere $\delta\cup d$ is embedded and bounds a 3-ball $b$ in $N_i$ (when the component of  $f^{-n}(S_{Q_i})$ containing $d$ is a 2-sphere, replace $d$ by the complementary disc if necessary). 

There are two occurrences: (a) $f^{n}(b)\subset Q_i$ and (b) $f^{n}(b)\subset f^{-1}(Q_i)$. We define $Q_i^{\prime}$ as $cl(Q_i\setminus f^{n}(b))$ in the case (a) and $Q_i\cup f^{n}(b)$ in the case (b).  The fact that $c$ is an innermost curve implies
$f(Q_i) \subset Q_i^{\prime}\subset Q_i$ in  case (a) and $Q_i\subset Q_i^{\prime}\subset f^{-1}(Q_i)$ in  case (b). In both cases there is a smooth approximation $\tilde Q_i$ of $Q_i^{\prime}$ 
such that $f(Q_i)\subset int~\tilde Q_i\subset Q_i$ if (a),  $Q_i\subset int~\tilde Q_i\subset f^{-1}(Q_i)$ if (b), and the 
number of intersection  curves in $S_{P_i}\cap (\bigcup\limits _ {n>0}f^ {-n} (\partial\tilde Q_i))$ is less than the cardinality of  $\mathcal C$. 

In  case (a), $\varphi_{_{f(Q_i)}}=\varphi_{_{Q_i}}f^{-1}:f(Q_i)\to\mathbb R$ is a dynamically ordered energy function which is constant on the boundary. 
Therefore $\tilde Q_i\setminus int~f(Q_i)$ is an 
$(f,S_{i})$-compressed product and, hence, 
due to proposition \ref{f,S}, there is a similar function on $\tilde Q_i$. Similarly in case (b), $\tilde Q_i$ is equipped with a dynamically ordered energy function which is constant on the boundary as 
$\tilde Q_i\setminus int~Q_i$ is an $(f,S_{i})$-compressed product.

We will repeat this process until getting a handle neighborhood $\hat Q_i$ of genus $g_i$  for the attractor $A_i$ such that 
$S_{P_i}\cap(\bigcup\limits_{n>0}f^{-n}
(\partial\hat Q_i))$ does not contain curves which are  homotopic zero in $S_{P_i}$. Thus we may assume that $S_{P_i}\cap(\bigcup\limits_{n>0}f^{-n}
({S_{Q_i}}))$ does not contain intersection curves which are homotopic to zero in $S_{P_i}$. 

We denote by $m$ the largest integer such that
$f^m(S_{P_i})\cap S_{Q_i}\neq\emptyset$. Let $F$ be a connected component of $f^m(S_{P_i})\cap G_i$.
We have $\partial F\subset\partial{S}_{Q_i}$.
Let us show that $F$ is incompressible in $G_i$. 
Indeed, if $\delta$ is a disc in $G_i$ with boundary $\gamma \subset F$ then $\gamma$ bounds 2-disk $\tilde\delta\subset f^{m}(S_{P_i})$ as  $f^{m}(S_{P_i})$ is incompressible surface in $V_i$. By assumption the components of $\partial F$ are not homotopic to zero in $f^{m}(S_{P_i})$, then $\partial F\cap\tilde\delta=\emptyset$ and, hence, $\tilde\delta\subset F$. 

Therefore, according to statement \ref{wald32} 
there is some surface $F_1\subset S_{Q_i}$
diffeomorphic to $F$, with $\partial
F=\partial F_1$, and $F\cup F_1$
bounds  a domain $\Delta$ in $G_i$ which, up to
smoothing of the boundary, is diffeomorphic to $F\times[0,1]$. We then define $\tilde Q_i$ as
$Q_i\cup \Delta$ up to smoothing. By the choice of $m$, $\tilde Q_i$ is $f$-compessed as $f(\Delta)\subset Q_i$. As $\tilde Q_i$ is obtained by an isotopy supported in a neighborhood of 
$\Delta$ from $Q_i$ then $\tilde Q_i\setminus int~Q_i$ is an $(f,S_{i})$-compressed product. Thus  we get a dynamically  ordered energy function on $\tilde Q_i$ 
with $\partial\tilde Q_i$ as a level set. Arguing recursively, we are reduced to case 1).
\end{demo}

\begin{lemm} \label{gi+1gi} Let $i\in\{k_0+1,\dots,k_1\}$, $M_i$ be strongly tight neighborhood of the attractor $A_i$,   $D_{i}=M_{i}\cap W^s_{\mathcal O_{i}}$ and  $N(D_{i})\subset M_{i}$ be a tubular neighborhood of $D_{i}$ such that $N(D_{i})\cap A_{i-1}=\emptyset$ and the set  $P_{i-1}=M_{i}\setminus int~N(D_{i})$ is $f$-compressed. Then $P_{i-1}$ is a handle  neighborhood of genus $g_{i-1}$ for the attractor $A_{i-1}$. 
\end{lemm}
\begin{demo} Similar to proposition \ref{ttt} it is proved the equality $\bigcap\limits_{k\geq 0}f^k(P_{i-1})=A_{i-1}$.  Thus $P_{i-1}$ is a trapping neighborhood of the attractor $A_{i-1}$ and, hence, the set $P_{i-1}$ consists of $c_{i-1}$ connected components. Each of them is handlebody, as it is obtained from $M_i$ removing $(r_i-r_{i-1})$ 1-handles, which are the set $N(D_i)$. As in proof of proposition \ref{ttt}, a sum $g_{_{P_{i-1}}}$ of genera $P_{i-1}$ is calculated by formula $g_{_{P_{i-1}}}=g_i-((r_{i}-r_{i-1})-(c_{i-1}-c_i))$. Then  $g_{_{P_{i-1}}}=c_i+r_i-s_i-((r_{i}-r_{i-1})-(c_{i-1}-c_i))=c_{i-1}+r_{i-1}-s_i$. As $s_{i-1}=s_i$ then $g_{_{P_{i-1}}}=g_{i-1}$. 
\end{demo}

\subsection{Global construction}

We divide a construction of the dynamically ordered energy function for $f:M\to M$ on steps. 

{\bf Step 1.} By induction on $i=1,\dots,k_1$ let us prove the existence of a dynamically ordered energy function  $\varphi_{_{M_i}}$ on $M_i$ of the attractor $A_i$ with level set $S_i$. 

For $i=1$ the attractor $A_1$ coincides with sink orbit $\mathcal O_1$ of the diffeomorphism $f$. According to statement \ref{loc} there is a neighborhood $U_{\mathcal{O}_{1}}\subset M$ of the orbit $\mathcal O_{1}$, equipped by an energy function $\varphi_
{\mathcal{O}_{1}}:U_{\mathcal{O}_{1}}\to\mathbb{R}$ for $f$ and such that   $\varphi_{\mathcal{O}_{1}}(\mathcal{O}_{1})=1$. Moreover, for each connected component $U_{\omega}$, $\omega\in \mathcal{O}_{1}$ of the set   $U_{\mathcal{O}_{1}}$ there are Morse coordinates  $(x_1,x_2,x_3)$ such that  
$\varphi_{\mathcal{O}_{1}}(x_1,x_2,x_3)=1+x_1^2+x_2^2+x_3^2$. Then there is a value $\varepsilon_1>0$ such that  set  $Q_1=\varphi^{-1}_{\mathcal{O}_{1}}(1+\varepsilon_1)$ consists of $f$-compressed union of $c_1$ 3-balls. Thus $Q_1$ is a handle neighborhood of genus $0$ for the attractor $A_1$. As $g_1=0$ then, according to lemma \ref{ww}, there is a dynamically ordered energy function $\varphi_{_{M_1}}$ on the neighborhood $M_1$ for the attractor $A_1$ with level set $S_1$. 

Let, by assumption of the induction, there is a dynamically ordered energy function $\varphi_{_{M_{i-1}}}$ on the neighborhood $M_{i-1}$ of the attractor $A_{i-1}$ with level set $S_{i-1}$. Let us construct the function $\varphi_{_{M_i}}$. There is two cases: a) $i\leq k_0$; b) $i>k_0$. 

In the case a) the neighborhood $M_i$ consists of a handle neighborhood of genus $0$ for the attractor $A_{i-1}$ (denote it $P_{i-1}$) and a trapping neighborhood of the orbit  $\mathcal O_{i}$, consisting from 3-balls (denote it $Q_i$). By assumption of the induction and lemma \ref{ww}   there is a dynamically ordered energy function $\varphi_{_{P_{i-1}}}$ on $P_{i-1}$ with level set  $\partial P_{i-1}$. Similar to case $i=1$ it is shown the existence of a dynamically ordered energy function $\varphi_{_{Q_i}}$ on $Q_i$ with level set  $\partial Q_i$.  The required function $\varphi_{_{M_i}}$ is formed from $\varphi_{_{P_{i-1}}}$ and $\varphi_{_{Q_i}}$. 

В случае b) we follow to scheme of proof from section  4.3 of paper \cite{GrLaPo}. 

According to statement \ref{loc}, the orbit  $\mathcal O_{i}$ has  a neighborhood $U_{\mathcal{O}_{i}}\subset M$ endowed with an energy function $\varphi_
{\mathcal{O}_{i}}:U_{\mathcal{O}_{i}}\to\mathbb{R}$ of $f$ with $\varphi_
{\mathcal{O}_{i}}(\mathcal{O}_{i})=i$. Moreover, each connected component $U_{\sigma}$, $\sigma\in \mathcal{O}_{i}$ of $U_{\mathcal{O}_{i}}$ is endowed with   Morse coordinates $(x_1,x_2,x_3)$ such that 
$\varphi_{\mathcal{O}_{i}}(x_1,x_2,x_3)=i+x_1^2+x_2^2+x_3^2$, the $x_1$-axis is
contained in the unstable manifold and the $(x_2,x_3)$-plane is contained in the stable manifold of $\sigma$.

It follows from properties of strongly tight neighborhood  $M_i$ and $\lambda$-lemma that there is a tubular neighborhood $N(D_{i})\subset M_{i}$ of  $D_{i}=M_i\cap W^s_{\mathcal O_i}$ such that $N(D_{i})\cap A_{i-1}=\emptyset$, set $P_{i-1}=M_{i}\setminus int~N(D_{i})$ is  $f$-compressed and surface $\partial P_{i-1}$ transversal intersects each connected component of the set $\varphi^{-1}_{\mathcal O_{i}}(i)\setminus\mathcal O_{i}$ at one closed curve. By lemma \ref{gi+1gi}, the set $P_{i-1}$ is a handle neighborhood of genus $g_{i-1}$ for the attractor $A_{i-1}$. By assumption of induction and lemma \ref{ww} there is a dynamically ordered energy function $\varphi_{_{P_{i-1}}}$ on $P_{i-1}$ with level set $\partial P_{i-1}$. 

\begin{figure}\epsfig
{file=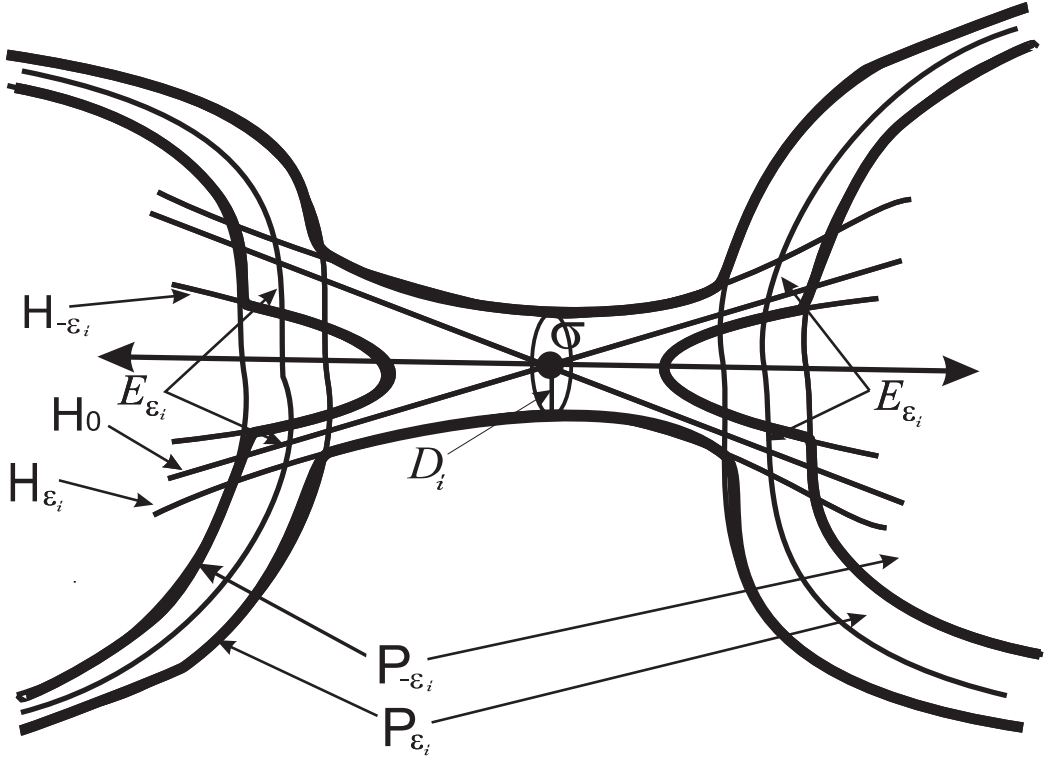, width=14. true cm, height=10 true
cm} \caption{Illustration to step 1}\label{nefig}
\end{figure}

For $\varepsilon_i\in(0,1),~t\in[-\varepsilon_i,\varepsilon_i]$ set $P_{t}=\varphi^{-1}_{{P_{i-1}}}([1,\varphi_{_{P_{i-1}}}(\partial P_{i-1})-\varepsilon_i+t])$, $H_t=\{x\in U_{\mathcal O_{i}}~:~\varphi_{\mathcal O_{i}}(x)\leq i+{t}\}$ and $E_{\varepsilon_i}=(P_{\varepsilon_i}\setminus int~P_{-\varepsilon_i})\cap(H_{\varepsilon_i}\setminus int~H_{-\varepsilon_i})$ (see figure  \ref{nefig}). Notice that $P_{\varepsilon_i}=P_{i-1}$ and, hence, $f(P_{\varepsilon_i})\subset int~P_{\varepsilon_i}$. As $\varphi_{\mathcal O_{i}}$ is a Lyapunov function for $f|_{U_{\mathcal O_{i}}}$ then $\varphi_{\mathcal O_{i}}(f^{-1}(\varphi^{-1}_{\mathcal O_{i}}(i)\setminus\mathcal O_{i}))>i$ and, hence, $(H_0\setminus\mathcal O_i)\subset int~f^{-1}(H_0\setminus\mathcal O_i)$. This and conditions of choice of  $N(D_{i})$  implies the existence of a value ${\varepsilon_i}$ with following properties:
\begin{itemize}
\item[(1)] $f(P_{\varepsilon_i})\subset int~P_{-\varepsilon_i}$;
\item[(2)] for each $t\in[-{\varepsilon_i},{\varepsilon_i}]$ surface  $\partial P_t$ transversal intersects each connected component of the set $\partial H_t\setminus D_{i}$ по at one closed curve;
\item[(3)] $f^{-1}(E_{\varepsilon_i})\cap H_{\varepsilon_i}=\emptyset$.
\end{itemize}

For $t\in[-\varepsilon_i,\varepsilon_i]$ set  $Q_t=P_t\cup H_t$. By the construction the set $Q_t,~t\neq 0$ is $f$-compressed. Moreover,  $Q_{-\varepsilon_i}$ after smoothing is a handle neighborhood of genus $g_{i-1}$ of the attractor $A_{i-1}$ and $Q_{\varepsilon_i}$ after smoothing is strongly tight  neighborhood of the attractor $A_{i}$. By assumption of induction and lemma \ref{ww} there is a dynamically ordered energy function $\varphi_{_{{Q}_{-\varepsilon_i}}}$ on ${Q}_{-\varepsilon_i}$, which is a constant on $\partial {Q}_{-\varepsilon_i}$. As $\varphi_{_{{Q}_{-\varepsilon_i}}}(A_{i-1})\leq i-1$ then, due to proposition \ref{f,S}, we can suppose that $\varphi_{_{{Q}_{-\varepsilon_i}}}({Q}_{-\varepsilon_i})=i-{\varepsilon_i}$.

Define function  
$\varphi_{_{{Q_{\varepsilon_i}}}}:Q_{\varepsilon_i}\to\mathbb R$ on the set  $Q_{\varepsilon_i}$ by formula:\\
$\varphi_{_{{Q_{\varepsilon_i}}}}(x)=\cases{\varphi_{_{{Q}_{-\varepsilon_i}}}(x),~x\in{Q}_{-\varepsilon_i};\cr i+t,~x\in Q_t.\cr}$ 
Let us check that  $\varphi_{_{{Q_{\varepsilon_i}}}}$ is a dynamically ordered energy function $f$, then the existence of required function $\varphi_{_{M_{i}}}:M_{i}\to\mathbb R$ will follow from lemma \ref{ww}. 

Represent the set ${Q}_{\varepsilon_i}$ as a union of  subsets with pairwise disjoint interiors: ${Q}_{\varepsilon_i}=A\cup B\cup C$, where  $A={Q}_{-\varepsilon_i}$, $B=P_{\varepsilon_i}\setminus{Q}_{-\varepsilon_i}$ and $C={Q}_{\varepsilon_i}\setminus (P_{\varepsilon_i}\cup{Q}_{-\varepsilon_i})$. By the construction $\varphi_{_{{Q_{\varepsilon_i}}}}|_{A}$ 
a dynamically ordered energy function for $f$,  $\varphi_{_{{Q_{\varepsilon_i}}}}(\partial{A})=i-\varepsilon_i$, the function $\varphi_{_{{Q_{\varepsilon_i}}}}|_{B}$ has no critical points and function $\varphi_{_{{Q_{\varepsilon_i}}}}|_{C}$ coincides with function $\varphi_{_{\mathcal O_i}}|_{C}$. 
Let us check decreasing property of  $\varphi_{_{{Q_{\varepsilon_i}}}}$ along trajectories of $f$.  

If $x\in A$ then $f(x)\in A$ and  $\varphi_{_{{Q_{\varepsilon_i}}}}(f(x))<\varphi_{_{{Q_{\varepsilon_i}}}}(x)$, as 
$\varphi_{_{{Q_{\varepsilon_i}}}}|_{A}$ is a Lyapunov function. If $x\in B$ then, due to condition (1) of choice of $\varepsilon_i$, $f(x)\in A$ and, hence, $\varphi_{_{{Q_{\varepsilon_i}}}}(x)>i-\varepsilon_i$,   $\varphi_{_{{Q_{\varepsilon_i}}}}(f(x))<i-\varepsilon_i$, therefor   $\varphi_{_{{Q_{\varepsilon_i}}}}(f(x))<\varphi_{_{{Q_{\varepsilon_i}}}}(x)$. If $x\in C$ then, due to condition (3) of choice of $\varepsilon_i$, either $f(x)\in A$ and decreasing is proved as for $x\in B$, or $f(x)\in C$ and decreasing follows from the fact that $\varphi_{_{{Q_{\varepsilon_i}}}}|_{C}$ is a Lyapunov function.

{\bf Step 2.} In this step we delive a construction similar to step 1 for diffeomorphism $f^{-1}$. For this aim we recall that dynamical numbering of the orbits $\mathcal O_1,\dots,\mathcal O_{k_f}$ of the diffeomorphism $f$ induces dynamical numbering of the orbits $\tilde{\mathcal O}_1,\dots,\tilde{\mathcal O}_{k_f}$ of the diffeomorphism  $f^{-1}$ following way: $\tilde{\mathcal O}_i=\mathcal O_{k_f-i}$. Denote by $\tilde A_i$ the attractors of the diffeomorphism $f^{-1}$, by $\tilde M_i$ thir neighborhood and by $\tilde g_i$ a number, defined by formula 
$\tilde g_i=\tilde c_i+\tilde r_i-\tilde s_{i}$, wher  $\tilde c_i$ the number of the connected components of the attractor $\tilde A_i$, $\tilde r_i$ the number of the saddle pointd and  $\tilde s_i$ the number of the sink points of the diffeomorphism $f^{-1}$, belonging to  $\tilde A_i$. 

Set $\tilde k_1=k_f-k_1$ and consider the attrator $\tilde A_{\tilde k_1}$ for the diffeomorphism $f^{-1}$ (which, recall, is a repeller for the diffeomorphism $f$). Similar to step 1 we construct a 
a dynamically ordered energy function  $\tilde\varphi_{{\tilde M_{\tilde k_1}}}$ for $f^{-1}$ on the neighborhood ${\tilde M_{\tilde k_1}}$ with level set $\tilde S_{\tilde k_1}=\partial {\tilde M_{\tilde k_1}}$. 

{\bf Step 3.} In this step we show that set $P_{k_1}=M\setminus int~\tilde M_{\tilde k_1}$ is a handle neighborhood of genus $g_{k_1}$ of the attractor $A_{k_1}$, this implies the existence of the required function $\varphi$. Indeed, by lemma \ref{ww}, the existence of a dynamically ordered energy function $\varphi_{_{M_{k_1}}}$ on the neighborhood $M_{k_1}$ of the attractor $A_{k_1}$ implies the existence of a dynamically ordered energy function for $\varphi_{_{P_{k_1}}}$ on  $P_{k_1}$ with level set $\partial P_{k_1}$. 
According to proposition \ref{f,S} the function $\varphi_{_{P_{k_1}}}$ we can construct such that  $\varphi_{_{P_{k_1}}}(\tilde S_{\tilde k_1})=k_f+1-\tilde\varphi_{{\tilde M_{\tilde k_1}}}(\tilde S_{\tilde k_1})$. 
As $\partial P_{k_1}=\tilde S_{\tilde k_1}$ then required function $\varphi$ is defined by formula   $\varphi(x)=\cases{ \varphi_{_{P_{k_1}}}(x),~x\in P_{k_1};\cr k_f+1-\tilde\varphi_{{\tilde M_{\tilde k_1}}}(x),~x\in{\tilde M_{\tilde k_1}}.\cr}$  

Thus, let us prove that the set $P_{k_1}=M\setminus int~\tilde M_{\tilde k_1}$ is a handle neighborhood of genus $g_{k_1}$ of the attractor $A_{k_1}$. Set 
$\tilde N_{\tilde k_1}=W^s_{\tilde A_{\tilde k_1}\cap\Omega_{f^{-1}}}$ and $\tilde V_{\tilde k_1}=\tilde N_{\tilde k_1}\setminus\tilde A_{\tilde k_1}$. Notice that the open sets $V_{k_1}$ and $\tilde V_{\tilde k_1}$ are coincide, as both are obtained from $M$ by removing of $A_{k_1}$ and $\tilde A_{\tilde k_1}$. It follows from proof of proposition \ref{coon} that each of next sets $A_{k_1},~\tilde A_{\tilde k_1},M_{k_1},~\tilde M_{\tilde k_1},~N_{k_1},~\tilde N_{\tilde k_1},~V_{k_1},~\tilde V_{\tilde k_1}$ is connected. Then $g_{k_1}=1+|\Omega_1|-|\Omega_0|$ and $\tilde g_{\tilde k_1}=1+|\Omega_2|-|\Omega_3|$. From statement \ref{g(f)} we get $g_{k_1}=\tilde g_{\tilde k_1}$. Thus the handle neighborhoods $M_{k_1}$ and  $\tilde M_{\tilde k_1}$ have the same genera and their boundaries $S_{k_1}$ and $\tilde S_{\tilde k_1}$ belong to the set $V_{k_1}$, which is diffeomorphic to $S_{k_1}\times\mathbb R$. 

Choose $n\in\mathbb N$ such that $f^n(M_{k_1})\subset int~P_{k_1}$. Then, according to ring hypothesis and corollary \ref{the3}, manifold $K=P_{k_1}\setminus int~f^n(M_{k_1})$ is diffeomorphic to  $S_{k_1}\times[0,1]$. By the construction $f^n(M_{k_1})$ is a handle neighborhood of genus $g_{k_1}$ of the attractor $A_{k_1}$ and $P_{k_1}=f^n(M_{k_1})\cup K$. This implies that $P_{k_1}$ also is  handle neighborhood of genus $g_{k_1}$ of the attractor $A_{k_1}$.

\section{Dynamically ordered energy function for diffeomorphisms on 3-sphere}

In this section $f:\mathbb S^3\to\mathbb S^3$ is  a Morse-Smale diffeomorphism without heteroclinic curve. 

{\bf Proof of theorem \ref{sph}} 

Let us prove that diffeomorphism $f$ possesses a 
dynamically ordered energy function if and only if all its one-dimensional attractors and repellers are  tightly embedded. 

\begin{demo} The necessity of conditions of the theorem follows from \ref{tt}, let us proof the sufficiency. 

Let $i=k_0+1,\dots,k_1$. Then $A_i$ is one-dimensional attractor, consisting of $c_i$ connected components, containing $r_i$ saddles,  $s_i$ sinks and for which a number $g_i$ can be calculated by formula $g_i=c_i+r_i-s_i$. Firstly prove that $g_i=0$ for each $i=k_0+1,\dots,k_1$. 

We start from $g_{k_1}$. According to proposition \ref{coon}, the attractor  $A_{k_1}$ is connected that is $m_{k_1}=1$ and, hence,   $g_{k_1}=1+|\Omega_1|-|\Omega_0|$. Due to statement  \ref{g(f)}, we have $g_{k_1}=\tilde g_{k_1}$, where  $\tilde g_{k_1}=1+|\Omega_2|-|\Omega_3|$. According to statement \ref{the}, $2+|\Omega_1\cup\Omega_2|-|\Omega_0\cup\Omega_3|=0$ for any Morse-Smale diffeomorphism without heteroclinic curves on $\mathbb S^3$. Thus $g_{k_1}+\tilde g_{k_1}=0$ and, hence, $g_{k_1}=\tilde g_{k_1}=0$. Further let us show that $g_{i}\leq g_{i+1}$ for each $i=k_0,\dots,k_1-1$. 

Indeed, $g_{i+1}-g_i=(c_{i+1}-c_i)+(r_{i+1}-r_i)-(s_{i+1}-s_i)$. At the same time  $(c_{i}-c_{i+1})\leq(r_{i+1}-r_i)$, $s_{i+1}=s_i$ and, hence, $g_{i+1}\geq g_i$.  

Thus,  $g_i=0$ for each $i=k_0+1,\dots,k_1$. Then $K_i=M_i\setminus int~f(M_i)$ is a union of 3-dimensional annulus $S^2\times [0,1]$. As  $M_i\setminus A_i=\bigcup\limits_{k\geq 0}f^k(M_i)$ then  $M_i\setminus A_i$ is diffeomorphic to  $\partial M_i\times(0,1]$. Thus, the attractor $A_i$ is strongly tight embedded. Similar fact has place for repellers. The, according to theorem \ref{iff}, $f$ possesses a dynamically ordered energy function.
\end{demo}

\end{document}